\documentclass[12pt, a4paper]{article}
\usepackage[T1]{fontenc}
\usepackage{amsmath,amsthm}
\usepackage{amssymb}
\usepackage{amscd}
\usepackage{euscript}
\usepackage[latin1]{inputenc}
\usepackage{amsfonts}
\usepackage[matrix,arrow,curve]{xy}

\theoremstyle{plain}
\newtheorem{theo}{Theorem}[section]
\newtheorem{lem}[theo]{Lemma}
\newtheorem{prop}[theo]{Proposition}

\theoremstyle{definition}

\newtheorem{theosd}[theo]{Theorem}
\newtheorem{rem}[theo]{Remark}

\newcommand{\Z}{\mathbb{Z}}

\newcommand{\FF}{\ensuremath{\mathbb{F}}}

\newcommand{\ZZl}{\ensuremath{\mathbb Z}/{\ell}}

\newcommand{\Hn}{\ensuremath{H_{\mathrm{nr}}}}
\newcommand{\PP}{\ensuremath{\mathcal{P}}}

\newcommand{\eqdef}{\ensuremath{\stackrel{\mathrm{def}}{=}}}

\title{On a local-global principle for $H^3$ of function fields of surfaces over a finite field.}

\author{Alena Pirutka }

\begin{document}
\maketitle

\begin{abstract}
Let $K$ be the function field of a smooth projective surface $S$ over a finite field $\FF$. In this article, following the work of Parimala and Suresh, we establish a local-global principle for the divisibility of elements in $H^3(K,\ZZl)$ by elements in $H^2(K,\ZZl)$, $l\neq car. K$.
\end{abstract}

\section{Introduction}

Let $K$ be a field of one of two following types:
\begin{itemize}
\item [(i)] $K$ is the function field of a smooth projective surface $S$ over a finite field $\FF$ of characteristic  $p$;
\item [(ii)] $K$ is the function field of  a regular (relative) curve, proper over a ring of integers of a $p$-adic field.
\end{itemize}

Recall that for a field $k$, the u-invariant $u(k)$ is defined as a maximal dimension of an anisotropic quadratic form over $k$ \cite{Ka}. Another arithmetical invariant which one can associate to $k$ is the period-index exponent: an integer $d$ (if it exists), such that for any element $\alpha$ of the Brauer group $Br\,k$ we have $\mathrm{ind}\alpha| (\mathrm{per}\alpha)^d$. Recall that $\mathrm{ind}\alpha$ and $\mathrm{per}\alpha$ have the same prime factors, so that for a fixed $\alpha\in Br\,k$ one can find an integer $d=d(\alpha)$ as above.

For $k=K$ of one of the types above, these invariants are now understood. More precisely, if $K$ is of type $(ii)$, Saltman \cite{S97, S98} showed that $\mathrm{ind}\alpha| (\mathrm{per}\alpha)^2$ for $(\mathrm{per}\alpha,p)=1$. For $K$ of type $(i)$, using the techniques of twisted sheaves and also some Saltman's results on the ramification of $\alpha$, Lieblich \cite{Li} established that $\mathrm{ind}\alpha| (\mathrm{per}\alpha)^2$ as well. For the $u$-invariant,  Parimala and Suresh \cite{PS09} established that $u(K)=8$  for $K$ of type $(ii)$ and  $p\neq 2$ (in the case $(i)$  one easily sees that $u(K)=8$ as well). This result has been also  obtained by different methods by Harbater, Hartmann and Krashen \cite{HHK} and Heath-Brown and Leep \cite{HB, Leep} (which also contains the case $p=2$). One of crucial steps in the proof of Parimala and Suresh is to establish a local-global principle for the divisibility of elements in $H^3(K,\ZZl)$ by symbols $\alpha\in H^2(K,\ZZl)$, $\ell\neq p$, which also uses Saltman's results \cite{S97, Sa07} on the classification of the ramification points of $\alpha$. In \cite{PS10}, they also apply such a local-global principle to establish the vanishing of the third unramified cohomology group of conic fibrations over $S$.

In this note, we follow the arguments of Parimala and Suresh and we establish the local-global principle in the general case for $K$ of type $(i)$. The main technical difficulty is that for $\alpha$ a symbol there is no so-called <<hot>> points in the classification of Saltman, which was also the case considered in \cite{Sa07, Sa08}. Our main result is the following:

\theosd\label{lg}{\it Let $K$ be the function field of a smooth projective surface $S$ over a finite field and let $\ell$ be a prime, $\ell\neq char(K)$. Assume that $K$ contains a primitive $\ell^{\mathrm{th}}$ root of unity. Let $\xi\in H^3(K,\ZZl)$ and $\alpha\in H^2(K,\ZZl)$ be such that the union $\mathrm{ram}_S(\xi)\cup \mathrm{ram}_S(\alpha)$ is a simple normal crossings divisor. Assume that for any point $x\in S^{(1)}$ there exists $f_x\in K_x^*$ in the field of fractions of the completion of the local ring  of $S$ at $x$, such that $\xi=\alpha\cup f_x$ in $H^3(K_x,\ZZl)$. Then there exists a function $f\in K^*$ such that $\xi=\alpha\cup f$ in $H^3(K,\ZZl)$.\\}

In section \ref{snot}, we first fix the notations and recall some elements of Saltman's approach on the ramification of $\alpha\in Br\,k$, then we give some additional properties of so-called <<hot>> points. This allows us to  deduce the  local-global principle \ref{lg} in section \ref{slg}.
\paragraph{Acknowledgements.} I would like to thank Parimala and Suresh for their comments and corrections on the first versions of this article.

\section{Classification of the ramification points, complements on 'hot' points}\label{snot}
\subsection {Notations and first properties}

\subsubsection{Residues and unramified cohomology}

Let  $A$ be a discrete valuation ring of rank one with fraction field $K$ and the residue field $\kappa$, let $\pi$ be a uniformizing parameter of $A$.  For any $i\geq 1$, $j\in \mathbb Z$ and $n$ an integer  invertible in $\kappa$ we have the residue maps in Galois cohomology $$H^i(K, \mathbb \mu_n^{\otimes j})\stackrel{\partial_{A}}{\to}H^{i-1}(\kappa, \mathbb \mu_n^{\otimes j-1}).$$
If $\partial_A(x)=0$ for $x\in H^i(K, \mathbb \mu_n^{\otimes j})$ we say that $x$ is {\bf unramified}. In this case we define the {\bf specialisation} $\bar x$ of $x$ as $\bar x=\partial_A(x\cup \pi)$.

If $A$ is a regular ring with fraction field $K$, an element $x\in H^i(K, \mathbb \mu_n^{\otimes j})$ is called {\bf unramified} on $A$ if it is unramified with respect to all discrete valuations corresponding to height one prime ideals in $A$.

Let $k$ be a field. For $L$ a function field over $k$, $n$ an integer  invertible in $k$,   $i\geq 1$ and $j\in \mathbb Z$ we denote
 \vspace{-0.1cm}
$$H_{nr}^i(L/k,\mathbb \mu_n^{\otimes j})=\bigcap_A\mathrm{Ker}[H^i(L, \mathbb \mu_n^{\otimes j})\stackrel{\partial_{A}}{\to}H^{i-1}(k_A, \mathbb \mu_n^{\otimes j-1})],$$
\vspace{-0.4cm}

\noindent where $A$ runs through all discrete valuation rings of rank one with $k\subset A$ and fraction field $L$. Here we denote by $k_A$ the residue field of $A$ and by $\partial_{A}$ the residue map.
If $X$ is an integral variety over $k$, we denote
$$\Hn^i(X, \mu_n^{\otimes j})\eqdef\Hn^i(k(X)/k, \mu_n^{\otimes j}),
$$
where $k(X)$ is the function field of $X$.
If $X$ is a smooth and projective variety, then we also have \begin{equation}\label{nrs}\Hn^i(X, \mu_n^{\otimes j})=\bigcap_{x\in X^{(1)}}\mathrm{Ker} \partial_x,\end{equation} where we denote $\partial_x\eqdef \partial_{\mathcal O_{X,x}}$ (see \cite{CT}).\\

We will use the following vanishing results for varieties over finite fields.

\prop\label{hnul}{Let $\FF$ be a finite field and let $\ell$ be a prime different from the characteristic of $\FF$.
\begin{itemize}
 \item [(i)] If $C/\FF$ is a smooth projective curve,  then $\Hn^{2}(C, \mu_{\ell})=0$.
 \item[(ii)] If $S/\FF$ is a smooth projective surface, then $\Hn^{3}(S, \mu_{\ell}^{\otimes 2})=0$.
\end{itemize}
}
\proof{For $(i)$ note that  $\Hn^{2}(C, \mu_{\ell})$ is the $\ell$-torsion subgroup of $Br\,C$ (see \cite{CT}), the group which is zero as $C$ is a smooth projective curve over a finite field \cite{Gr}. The statement $(ii)$ is established in \cite{CTSS} p.790. \qed}
\rem{A part of the Kato conjecture states that $\Hn^{d+1}(X, \mu_{\ell}^{\otimes d})=0$ for $X$ a smooth  projective variety of dimension $d$, defined over a finite field $\FF$. This conjecture has been recently established by Kerz and Saito \cite{KS}, using also the arguments of Jannsen \cite{JS}. \\}

\subsubsection{Function fields of surfaces over a finite field}
In the rest of  this section we use the following notations :
 $\FF$ is a finite field of characteristic $p$, $\ell$ is a prime different from $p$,
$S$ is a smooth projective surface over $\FF$ and
$K$ is the function field of $S$.
If $P$ is a point of $S$, $\hat A_P$ denotes the completion of the local ring $A_P\eqdef\mathcal O_{S,P}$ at $P$ at its maximal ideal, $K_P$ is the field of fractions of $\hat A_P$.\\
\vspace{-0.4cm}
\begin{quote}
 {\it We assume that $K$ contains $\ell^{\mathrm{th}}$ roots of unity.}
\end{quote}
\noindent We denote by $\alpha$ a (fixed for what follows) element of  $ H^2(K,\ZZl)$.
For any $x\in S^{(1)}$ a codimension one point of $S$ we have the residue maps in Galois cohomology, as in the previous section:
$$\partial_x: H^i(K,\ZZl)\to H^{i-1}(\kappa(x),\ZZl),\; i\geq 1$$ where $\kappa(x)$ is the residue field of $x$;
the set of points $x\in S^{(1)}$ such that $\partial_x(\alpha)$ is non zero is finite.
We define the {\bf ramification divisor} of $\alpha$ as $\mathrm{ram}_S\alpha=\sum\limits_{x\in S^{(1)}, \partial_x(\alpha)\neq 0} x$.

We write $$\mathrm{ram}_S\alpha=\sum\limits_{i=1}^m C_i$$ where $C_i\subset S, i=1\ldots m$ are integral curves and we denote
$u_i=\partial_{C_i}(\alpha)$, $u_i\in H^1(\kappa(C_i),\ZZl)=\kappa(C_i)^*/\kappa(C_i)^{*\ell}$.

If $P\in S$ is a closed point, we will say that $\alpha$ is {\bf unramified at $P$} if $P\notin \mathrm{ram}_S\alpha$.

\lem\label{locnul}{Let $P\in S$ be a closed point. If $\alpha$ is unramified at $P$, then $\alpha$ is trivial in $H^2(K_P,\ZZl)$. }
\proof{As $\alpha$ is unramified at $P$, we have that $\alpha$ comes from $H^2_{\acute{e}t}(A_P, \ZZl)$ (cf. \cite[\S 3.6, \S 3.8]{CT}). It is then trivial in $H^2(K_P,\ZZl)$, as  the group $H^2_{\acute{e}t}(\hat A_P, \ZZl)=H^2(\kappa(P), \ZZl)$ is zero as $\kappa(P)$ is a finite field. \qed \\}

{\it In the rest of this section  we assume that $\mathrm{ram}_S\alpha$ is a simple normal crossings divisor.}

\subsection{Classification of points}
Let $P\in S$ be a closed point. We recall the classification of \cite{Sa07} with respect to the divisor $\mathrm{ram}_S\alpha:$
\begin{enumerate}
\item if $P\notin \mathrm{ram}_S\alpha$ then it is called a {\bf neutral} point;
\item if $P$ is on only one irreducible curve in the ramification divisor, then it is called a {\bf curve} point.
\item As the divisor $\mathrm{ram}_S\alpha$ is a simple normal crossings divisor, the remaining case is when $P$ is only on two curves of the ramification divisor, in this case it is called a {\bf nodal} point. Assume that $P$ is on the curves $C_i$ and $C_j$ for some $i\neq j$. Recall that we denote $u_i=\partial_{C_i}(\alpha).$
    By the reciprocity law (see \cite{Ka86}), $\partial_P(u_i)=-\partial_P(u_j)$. Then the following cases may occur:
\begin{enumerate}
\item $u_i$ and $u_j$ are ramified at $P$. Then $P$ is called a {\bf cold} point;
\item $u_i$ and $u_j$ are unramified at $P$. Denote $u_i(P), u_j(P)\in H^1(\kappa(P), \ZZl)$ the specialisations of $u_i$ (resp. $u_j$) at $P$. \begin{enumerate} \item If $u_i(P)$ and $u_j(P)$ are trivial, $P$ is called a {\bf cool} point;
    \item if $u_i(P)$ and $u_j(P)$ are non trivial and generate the same subgroup of  $H^1(\kappa(P), \ZZl)$, then $P$ is called a {\bf chilly} point;
      \item if $u_i(P)$ and $u_j(P)$  do not generate the same subgroup of  $H^1(\kappa(P), \ZZl)$, then $P$ is called a {\bf hot} point. Since $\kappa(P)$ is a finite field , we have $H^1(\kappa(P), \ZZl)\simeq \ZZl$. We then get that for a hot point one of the specialisations $u_i(P)$ and $u_j(P)$, say $u_j(P)$, should be trivial. Then we will call $P$ a {\bf hot non neutral point} on $C_i$ and {\bf hot neutral point} on $C_j$.
    \end{enumerate}
\end{enumerate}
\end{enumerate}

Consider a graph whose vertices are  curves in the ramification locus and whose edges correspond to chilly or hot points. We will call  a {\bf hot chilly circuit} a loop in this graph or a connected component with more than one hot point on curves in this component.

\subsection{Local description}

\lem\label{locdescr}{Let $P\in S$ be a closed point.\begin{itemize} \item[(i)] If $P\in C_i$ is a curve point with  $\pi_i$ a  local parameter of $C_i$ in $K$, then in $H^2(K_P,\ZZl)$ we have $\alpha=u\cup \pi_i$ for some unit $u\in A_P$.
\item[(ii)]
If $P\in C_i\cap C_j$ is a nodal point with  $\pi_i, \pi_j$  local parameters of $C_i$ and $C_j$ in $K$, then in $H^2(K_P,\ZZl)$ we have
\begin{enumerate}
\item $\alpha=0$ if $P$ is a cool point;
\item $\alpha=u\pi_j\cup v\pi_i^s$ for some units $u, v\in A_P$, $1\leq s\leq \ell-1$, if $P$ is a cold point;
\item $\alpha=u\cup \pi_i\pi_j^s$ for some unit $u\in A_P$, $1\leq s\leq \ell-1$, if $P$ is a chilly point;
\item $\alpha=u\cup \pi_i$, for some unit $u\in A_P$, if $P$ is a hot point, non neutral on $C_i$.
\end{enumerate}
\end{itemize}}
\proof{For (i) we  write  (see \cite[Proposition 2.1]{S97})
$$\alpha=\alpha'+(u,\pi_i)$$ where $\alpha'$ is unramified at $P$ and $u$ is a unit in $A_P$.  By lemma \ref{locnul},  $\alpha'$ is zero in $H^2(K_P,\ZZl)$, so that we get $(i)$.  For $(ii)$, all the statements but the last one are in \cite[Lemma 1.3]{PS10}. To see the last one, we  write  (see \cite[Proposition 2.1]{S97})
$$\alpha=\alpha'+(u,\pi_i)+(v, \pi_j)$$ where $\alpha'$ is unramified at $P$, $u,v$ are units in $A_P$
and the image of $v$ in $\kappa(C_j)$ is $\partial_{C_j}(\alpha)$. As above,  $\alpha'$ is zero  by lemma \ref{locnul}.  The element $v$  is also zero in $\hat A_P$ because $v(P)$ is zero in $H^1(\kappa(P), \Z/\ell)$  by the definition of a hot point, so that $v\in H^1(\hat A_P,\Z/\ell)=H^1(\kappa(P), \Z/\ell)$ is zero as well.\qed \\ }

Let $C_i$ be a curve in the ramification divisor $\mathrm{ram}_S\alpha$ and let $\pi_i$ be a prime defining $C_i$ in $K$. We define a {\bf residual class} $\beta_i$ of $\alpha$ at $C_i$ as
$\beta_i=\partial_{C_i}(\alpha\cup \pi_i)$ (see \cite[Remark 2.6]{PS10}, \cite[p.820]{Sa07}).\\ %This definition depends on a choice of $\pi_i$, but any residual class $\beta_i$ satisfies the following universal property: (see \cite[Lemma 2.5]{PS10}).
\prop\label{resnhp}{If $C_i\in \mathrm{ram}_S\alpha$ contains no hot points, the element $\beta_{i}$ is trivial over $L_i=\kappa(C_i)(\sqrt[\ell]u_i)$.}
\proof{
 Using proposition \ref{hnul}, it is sufficient to see that $\beta_i$ is unramified over $L_i$. Let $v$ be a discrete valuation on $L_i$, let $v'$ be the induced valuation on $\kappa(C_i)$ and let $e$ be the valuation of a uniformizing parameter of $\kappa(C_i)$ in $L_i$. We have the following commutative diagram (cf. \cite[Proposition 3.3.1]{CT}):
$$\xymatrix{
H^2(L_i,\ZZl)\ar[r]^{\partial_{v}} & H^1(\kappa(v),\ZZl)\\
H^2(\kappa(C_i),\ZZl)\ar[r]^{\partial_{v'}}\ar[u]^{res} & H^1(\kappa(v'),\ZZl)\ar[u]^{e\cdot res}.\\
}
$$

We may assume that the valuation $v'$ corresponds to a point $P\in C_i$. We have the following cases:
\begin{enumerate}
\item if $P$ is a curve point, in $\hat A_P$ we have $\alpha=u\cup \pi_i$ for some unit $u$ by the local description \ref{locdescr}, and so $\partial_{v'}(\beta_i)=\partial_P\partial_{C_i}(u\cup \pi_i\cup \pi_i)=0$.
\item if $P$ is a cool point, $\alpha$ is trivial in  $\hat A_P$ by the local description, hence  $\partial_{v'}(\beta_i)=0$.
\item if $P$ is a cold point, the extension $L_i/\kappa(C_i)$ is ramified at $P$, hence $l|e$. As $\kappa(v')$ is a finite field, $\partial_{v'}(\beta_i)$ becomes trivial in $\kappa(v)$. % we have in  $A_P$: $\alpha=(u\pi_j\cup v\pi_i^s)+\alpha'$ where $\alpha'$ is unramified at $P$ and $u,v$ are units (see \cite{Sa97} 2.1).  Then $u_i=\partial_{C_i}(\alpha)$ equals to the image of $u\pi_j$ in $\kappa(C_i)$. Hence $d=\partial_P(u\pi_j\cup v)$ is trivial over $L_i$.
\item  if $P$ is a chilly point, from the local description we deduce that $\partial_{v'}(\beta_i)=u_j(P)$ which equals to $u_i(P)^s$ for some $1\leq s\leq l-1$ by the definition of a chilly point, hence trivial over $\kappa(v)$.
\end{enumerate}
Thus $\beta_i$ is unramified over $L_i$ and hence trivial.
\qed\\}

\section{The local-global principle}\label{slg}
In this section we prove the local-global principle \ref{lg}.
\subsection{Additional notations}
Let $\xi\in H^3(K,\ZZl)$ be as in the theorem \ref{lg}. As the divisor $\mathrm{ram}_S(\xi)\cup \mathrm{ram}_S(\alpha)$ is a simple normal crossings divisor,  we can write $$\mathrm{ram}_S(\xi)\cup \mathrm{ram}_S(\alpha)=\sum\limits_{i=1}^n C_i$$ where $n\geq m$ and $C_i$ are integral curves.

We denote $\mathcal C=\{C_1,\ldots C_m\}$ and
$\mathcal T=\{C_1,\ldots C_m, C_{m+1}, \ldots C_{n}\}.$

Let $\PP$ a finite set of closed points on $S$, consisting of all the points of  intersections $C_i\cap C_j$
and at least one point from each component $C_i$, let
$B$ be the semilocal ring of $S$ at $\PP$ %$\{C_1,\ldots C_m, C_{m+1}, \ldots C_{n}\}$
and let $\pi_i$ be a (fixed from the beginning) prime defining $C_i$ in $B$.

Recall that we denote
$u_i=\partial_{C_i}(\alpha)\in \kappa(C_i)$. We put as before
$$\beta_i=\partial_{C_i}(\alpha\cup \pi_i)\in H^2(\kappa(C_i),\Z/\ell)$$
and
\begin{equation}\label{tj}f_i=f_{C_i}=\pi_i^{t_i}f_i'\in K_{C_i},\end{equation} where $v_{C_i}(f_i')=0.$

\subsection{Plan of the proof}

Vanishing of the third unramified cohomology of $S$ from proposition \ref{hnul} implies that to establish theorem \ref{lg} it is sufficient to prove that there exists $f\in K$
 such that  for any  $
 x\in S^{(1)}$, $\partial_{x}(\xi)=\partial_{x}(\alpha\cup f)$ in the residue field $\kappa(x)$. This is done essentially in two steps.
 \begin{enumerate}
 \item First we get the condition above for $x=C_i$; to do this we first work locally around the double points $P$.
 \item Next we adjust the element $f$ from the previous step to get the condition for other curves in the support of $\mathrm{div}(f)$.
 \end{enumerate}

 \subsection{Analyzing curves and points in the ramification divisor}

 \subsubsection{Local description at nodal points}
\prop\label{locan}{ Let $P\in C_i\cap C_j$ be a nodal point. Then
\begin{enumerate}
\item If $P$ is a chilly point, then for any $0\leq r_i\leq \ell-1$ there exists $0\leq r_j\leq \ell-1$ such that
$\xi-\alpha\cup \pi_i^{r_i}\pi_j^{r_j}$ is unramified on $\hat A_P$;
\item If $P$ is a hot point, non neutral on $C_i$, then
$\xi-\alpha\cup \pi_i^{r_i}\pi_j^{t_j}$ is unramified on $\hat A_P$ for any $0\leq r_i\leq \ell-1$\footnote{Moreover, one can see that for $r_j\neq t_j$ mod $\ell$, $\xi-\alpha\cup \pi_i^{r_i}\pi_j^{r_j}$ is ramified on $\hat A_P$}, where $t_j$ is defined in (\ref{tj}).\\
\end{enumerate}}
\proof{
The first part is in \cite[Lemma 3.1]{PS10}. To see the second, let us first write $t_j'=v_{\pi_j}(f_i)=v_{\pi_j}(f_i')$.

By the reciprocity law (cf. \cite{Ka86}), we have $$\partial_P(\partial_{C_i}(\xi))=-\partial_P(\partial_{C_j}(\xi)),$$ and so $$\partial_P(\partial_{C_i}(\alpha\cup f_i))=-\partial_P(\partial_{C_j}(\alpha\cup f_j)).$$
In $\hat A_P$ we have $\alpha=u\cup \pi_i$ by the local description.
This gives
$$\partial_P(\partial_{C_i}(u\cup \pi_i\cup \pi_i^{t_i}f_i'))=-\partial_P(\partial_{C_j}(u\cup \pi_i\cup \pi_j^{t_j}f_j')).$$

The left side is $$\partial_P(u(C_i)\cup f_i'(C_i)^{-1})=-(t_j'+\epsilon)u(P),$$
where we  denote $u(C_i),  f_i'(C_i)$ the images of $u, f_i'$ in $\kappa(C_i)$ and $\epsilon=v_P({f_i'(C_i)}/{\pi_j^{t_j'}}(C_i))$.

The right side is $$-t_j\partial_P(u(C_j)\cup \pi_i(C_j))=-t_ju(P),$$ where we  denote $u(C_j),  \pi_i(C_j)$ the images of $u, \pi_i$ in $\kappa(C_j)$.
We then get $t_j=t_j'+\epsilon$ as $u(P)$ is not an $l^{\mathrm{th}}$ power. We then can write $f_i=\pi_i^{t_i}\pi_j^{t_j-\epsilon}f_i''$, where $v_{P}(f_i''(C_i))=\epsilon$.\\

{\it Claim:} $\partial_{C_i}(\xi-\alpha\cup \pi_i^{r_i}\pi_j^{t_j})=0$ in the completion $\kappa(C_i)_P$.\\
In fact, in $\kappa(C_i)_P$, we  have
\begin{multline*}
\partial_{C_i}(\xi)=\partial_{C_i}(u\cup \pi_i\cup \pi_i^{t_i}\pi_j^{t_j-\epsilon}f_i'')=
\partial_{C_i}(u\cup \pi_i\cup \pi_i^{t_i}\pi_j^{t_j})+\partial_{C_i}(u\cup\pi_i\cup \pi_j^{-\epsilon} f_i'')=\\=
\partial_{C_i}(\alpha\cup \pi_i^{r_i}\pi_j^{t_j})-u(C_i)\cup \pi_j(C_i)^{-\epsilon}f_i''(C_i)
\end{multline*}
%u(C_i)\cup \pi_j^{-t_j}f_i''(C_i)^{-1}.$$
(for the last equality we observe that $\partial_{C_i}(u\cup \pi_i\cup \pi_i$) iz zero in $\kappa(C_i)_P$)
and $u(C_i)\cup \pi_j(C_i)^{-\epsilon}f_i''(C_i)$ is unramified at $P$ and, hence, zero in   $\kappa(C_i)_P$.  %$$\partial_{C_i}(\alpha\cup \pi_i^{r_i}\pi_j^{t_j})=\partial_(C_i)(u\cup \pi_i\cup \pi_i^{r_i}\pi_j^{t_j})=u(C_i)\cup \pi_j^{-t_j}.$$

From the claim we deduce that $\xi-\alpha\cup \pi_i^{r_i}\pi_j^{t_j}$ is ramified on $\hat A_P$ at most at $\pi_j$; using reciprocity (cf.\cite[Lemma 1.2]{PS10}), it is then unramified. \qed \\}

Using the proposition above, we would like to associate an integer $r_i$ to each curve $C_i$ globally, and not only locally in each double point. To do this, we will need to avoid the hot chilly circuits.

\lem\label{bup}{Let $P\in S$ be a chilly or a hot nodal point. Let $S'\to S$ be the blowing-up of $S$ at $P$. Then for any point $x\in S'^{(1)}$ there exists $f_x\in K_x^*$  such that $\xi=\alpha\cup f_x$ in $H^3(K_x,\ZZl)$.}
\proof{We need only to consider the case when $x=E$ is an exceptional divisor of the blowing-up. If $P$ is a chilly point, then the statement is established in \cite[Theorem 3.4]{PS10}. Assume that $P$ is a hot point, which is non neutral on $C_i$. Denote again $C_i,C_j$ the strict transforms of $C_i$ and $C_j$, let $P_i=C_i\cap E$ and $P_j=C_j\cap E$. Using the local description, we easily check that $P_j$ is a hot point, non neutral on $E$ and that $P_i$ is a chilly point. Then, taking $f_E=\pi_i^{r_i}\pi_j^{t_j}$ we see that $\xi-\alpha\cup f_E$  is unramified on $K_E$ (for any $r_i$). In fact, we have an injection of local rings $A_P\to \mathcal O_{S', E}$, so that $\hat A_P$ is contained on the completion of $\mathcal O_{S', E}$ at its maximal ideal. Since $\xi-\alpha\cup f_E$ is unramifed on $\hat A_P$ by proposition \ref{locan}, it is unramified at $E$. \qed\\}

\lem\label{nocircuits} {There exists a smooth and projective surface $S'\to S$, obtained as a blowing-up of $S$ in some hot and chilly points, such that on $S'$ there is no hot chilly circuits for $\alpha$ and such that we still have that for any point $x\in S'^{(1)}$ there exists $f_x\in K_x^*$  with $\xi=\alpha\cup f_x$ in $H^3(K_x,\ZZl)$. }
\proof{
In fact, suppose that we have a hot chilly circuit and let $P\in C_i\cap C_j$ be a hot or a chilly point on a curve in this circuit. If $P$ is a hot point, which is non neutral on $C_i$, let us consider the blowing-up of $P$. Let $E$ be the exceptional divisor of the blowing up at $P$, and denote again $C_i,C_j$ the strict transforms of $C_i$ and $C_j$, let $P_i=C_i\cap E$ and $P_j=C_j\cap E$. As in the previous lemma,   $P_j$ is a hot point, non neutral on $E$ and  $P_i$ is a chilly point.
We may then suppose that there is a chilly point on curves in the circuit. By a construction of \cite[2.9]{Sa07} we may break the circuit  by successif blowing-ups of this point, and (some of) the double points above. Thus by successif blowing-ups of hot and chilly points we can obtain a surface $S'$ with no hot chilly circuits. We may find a corresponding function $f_x$ for each exceptional divisor by the previous lemma.\\

{\it Up to replacing $S$ by $S'$ as in the previous lemma, from now on, we assume that there is no hot chilly circuits.\\}

Now we can choose constants $r_i$ globally.

\lem\label{rglob}{ Assume that there is no hot chilly circuits. Then we may  associate to each curve  $C_i\in \mathcal C$  an integer $0\leq r_i\leq \ell-1$  such that $\xi-\alpha\cup \pi_i^{r_i}\pi_j^{r_j}$ is unramified on $\hat A_P$ for any nodal point $P$. If there is a hot point on $C_i$ then $r_i=t_i$.}\proof{Since there is no hot chilly circuits,  we may successively choose the integers $r_i$ as in the statement using lemma \ref{locan}. \qed \\}

\subsubsection{First choice of $f$}
For $1\leq i\leq m$, let $r_i$ be as in lemma \ref{rglob}. For $i\geq m+1$, we put $r_i=v_{C_i}(f_{C_i})$.
\prop\label{1choice}{Assume that there is no hot chilly circuits. Then there exists a function
 $f\in K^*$ such that $$\mathrm{div}(f)=\sum_{i=1}^n r_iC_i+F,$$ where $F$ does not pass through any point in $\mathcal P$, and such that $\xi-\alpha\cup f$ is unramified on $\hat A_P$ for all $P\in \mathcal P$.}\proof{By \cite[Lemma 3.2]{PS10}, for any  non hot point $P\in C_i\cap C_j$  for $C_i,C_j\in\mathcal T$, one can associate $w_P\in A_P$,  such   that $\xi-\alpha\cup w_P\pi_i^{r_i}\pi_j^{r_j}$ is unramified at $\hat A_P$, where the $r_i$'s are as above. We put $w_P=1$ if $P$ is a hot point. Then $\xi-\alpha\cup w_P\pi_i^{r_i}\pi_j^{r_j}$ is unramified at $\hat A_P$ by lemma \ref{rglob}. For any other point $P\in \PP$ we put $w_P=1$. Next we choose the function $f$ satisfying the conditions of the proposition by the same construction as in  \cite[Lemma 3.3]{PS10}.  We recall shortly this construction. First consider $g=\prod\limits_{i=1}^r \pi_i^{r_i}$. Then we can choose an element $u_P\in A_P$ for any $P\in \PP$ as follows. If $P\notin C_i$ for all $i$, we put $u_P=1$. If $P$ is on only one $C_i$ we put $u_P=g\pi_i^{-r_i}$ and if $P\in C_i\cap C_j$ we put $u_P=g\pi_i^{-r_i}\pi_j^{-r_j}$. Now let $w\in B$ be such that $w(P)=w_P/u_P$ for all $P\in \PP$. Then one checks that $f=wg$ satisfies the conditions of the proposition.   \qed\\}

\subsubsection{A better choice of $f$}
To choose an adjusting term  we will need to make a more precise choice of the function $f$:

\prop\label{2choice}{After possibly blowing up the surface $S$, there exists an element $f\in K^*$ such that
\begin{itemize}
\item[(i)] $\mathrm{div}(f)=\sum_{i=1}^n r_iC_i+E;$
\item[(ii)] $\xi-\alpha\cup f$ is unramified at $C_i$ for any $C_i\in \mathcal T$;
\item[(iii)] $E$ does not pass through any point of $\mathcal P$;
\item[(iv)] $(E\cdot C_i)_P$ is a multiple of $\ell$ at any intersection point of $E$ and $C_i$, $i=1,\ldots n$.
\end{itemize}
}
\proof{Let us consider $f$ as in  proposition \ref{1choice}. First we claim that for $i\leq m$ the residue $$\gamma_i:=\partial_{C_i}(\xi-\alpha\cup f)$$ may be written as $u_i\cup b_i$ for some constant $b_i$.
In fact,  it is sufficient to prove that  $\gamma_i$ is trivial over $\kappa(C_i)(\sqrt[\ell]u_i)$.
We  have $\partial_{C_i}(\xi)=\partial_{C_i}(\alpha\cup f_i)=t_i\beta_i+u_i\cup g_i$ for some $g_i$
and $\partial_{C_i}(\alpha\cup f)=r_i\beta_i+u_i\cup h_i$ for some $h_i$. If $C_i$ has hot points, then $r_i=t_i$ by construction and the claim is clear.
If $C_i$ has no hot points, then $\beta_i$ is trivial over $\kappa(C_i)(\sqrt[\ell]u_i)$ by proposition \ref{resnhp}, hence  so is $\gamma_i$.

Note that for $i>m$ the residu $\gamma_i:=\partial_{C_i}(\xi-\alpha\cup f)$ is trivial by the choice of $f$.

We proceed next exactly as in \cite[Theorem 3.4]{PS10}. We have
\begin{enumerate}
\item $\xi-\alpha\cup f$ is unramified at $\hat A_P$ for all $P\in \mathcal P$ by the choice of $f$, so $\gamma_i$ is unramified at all points $P\in \mathcal P\cap C_i$. Hence $b_i$ is a norm from $\kappa(C_i)_P(\sqrt[\ell]u_i)$.
\item By weak approximation, we then can find $a_i\in\kappa(C_i)$ which is a norm from $\kappa(C_i)(\sqrt[\ell]u_i)$ and such that $a_ib_i(P)=1$ for all $P\in \mathcal P\cap C_i$. As $u_i\cup b_i=u_i\cup a_ib_i$ we may assume that $b_i(P)=1$ for all $P\in \mathcal P\cap C_i$.
\item   We take $b\in B$ a unit such that $b_i$ is the image of $b$ in $B/\pi_i$.
\item Changing $f$ by $bf$, the conditions $(i)-(iii)$ are now satisfied. By \cite[Theorem 3.4. p.15]{PS10}   $(E\cdot C_i)_P$ is a multiple of $\ell$, after possibly some blowing ups of the points in the intersection of $E$ and $C_i$. \qed \\
\end{enumerate}}

 \subsection{Adjusting term for other curves}
 Let $f$ and $E$ be as in proposition \ref{2choice}. Let $\mathcal P'$ be  a finite set of closed points on $S$, consisting of $\mathcal P$ and all intersection points of $C_i$ and $E$,
and at least one point from each component of $E$ and at least one non hot point from each $C_i$.  \\

\prop\label{ct}{There exists $u\in K^*$ and $x\in K^*$ a norm from $K(\sqrt[\ell]u)$, such that
\begin{itemize}
\item[(i)] $\mathrm{div}(u)=-E+E'+\ell U$, where $E'$ does not pass through any point in $\mathcal P'\setminus (\bigcup\limits_{i\neq j} C_i\cap C_j)$;
\item[(ii)] the image  of $u$ in $\kappa(C_i)^*/\kappa(C_i)^{*\ell}$ equals to $u_i$;\\
\item[(iii)]  $\mathrm{div}(x)=-E+E''+\ell U'$ where $E''$ does not pass through any point in $\mathcal P'$ and any intersection point $C_i\cap E'$;
\item[(iv)] if $D$ is an irreducible curve in the $Supp(E'')$, then the specialization of $\alpha$ at $D$ is unramified at every discrete valuation of $\kappa(D)$ centered on a closed point of $D$.
\end{itemize}
}\proof{Let $\PP''=\mathcal P'\setminus (\bigcup\limits_{i\neq j} C_i\cap C_j)$ and let $B'$ be the semilocal ring of $S$ at $\PP''$. By \cite[Proposition 0.3]{Sa08}, one can find $v\in B'$ such that $v=u_i\text{ mod }\pi_i$. By \cite[Lemma 2.1]{PS10}, under the assumption that $(E\cdot C_i)_P$ is a multiple of $\ell$ at any intersection point of $E$ and $C_i$, $i=1,\ldots n$, one can find an element $z\in K^*$ such that:
\begin{enumerate}
\item $z$ is a unit at $C_i$ and maps to an $\ell^{\mathrm{th}}$ power in $\kappa(C_i)$ for all $i$;
\item $\mathrm{div}(z)=-E+Z_1+\ell Z_2$ where the support of $Z_1$ does not contain any point in $\PP$.
\end{enumerate}
Then the element $u=vz$ satisfies the conditions $(i)-(ii)$ of the proposition. By \cite[Lemma 2.3]{PS10} we can construct $x$ satisfiying $(iii)$. One then verifies $(iv)$ as in \cite[Lemma 2.4]{PS10}.  \qed \\}

 \subsubsection{End of the proof of theorem \ref{lg}}

Let $x$ be as in  proposition \ref{ct}. We will now show that $\xi-\alpha\cup (fx)$ is unramified at any codimension one point $D\in S$. We have $$\mathrm{div}(fx)=\sum_{i=1}^n r_iC_i+E''+\ell U'.$$
\begin{enumerate}
\item If $D\in \mathcal T$, then $\xi-\alpha\cup f$ is unramified at $D$ by the previous section.
If $D=C_i$ with $i\leq m$, then $\alpha=\alpha'+(u,\pi_i)$, where $\alpha'$ is unramified on $C_i$, by the choice of $u$. Hence $\partial_{D}(\alpha\cup x)=\partial_{C_i}(u\cup \pi_i\cup x)=0$ as $x$ is a norm from $K(\sqrt[\ell]u)$.
If $D=C_i$ with $i>m$, then $\alpha$ and $x$ are unramified at $D$, then so is $\alpha\cup x$.
\item If $D\notin \mathcal T\cup Supp(fx)$ or $D\in Supp(U')$ then $\xi-\alpha\cup (fx)$ is unramified at $D$.
\item Assume now that $D\in Supp(fx)\setminus (\mathcal T\cup U')$, i.e. $D\in Supp (E'')$. We have that $\xi$ and $\alpha\cup f$ are unramified at $D$.
So it is sufficient to show that $\partial_D(\alpha\cup x)=0$ which follows from  proposition \ref{ct}(iv). \qed \\
\end{enumerate}

\end{document}